\documentclass[12pt,twoside]{article}
\usepackage[english]{babel}     
\usepackage[cp850]{inputenc}    
\usepackage[T1]{fontenc}        
\usepackage{amsmath}
\usepackage{amssymb}
\usepackage{a4wide}             
\usepackage{tu-preprint}
\newtheorem{lem}{Lemma}
\newtheorem{theo}{Theorem}
\newtheorem{corr}{Corollary}
\newtheorem{defi}{Definition}

\newcommand{\Int}{\mathrm{int}}

\def\eps{\varepsilon}

\def\phi{\varphi}
\def\rho{\varrho}

\def\prf{\textit{Proof.\ }}
\def\qed{\hfill \quad\rule{1.5ex}{1.5ex}}

\def\lra{\longrightarrow}
\newcommand{\konv}{\mathop{\lra}\limits}

\newcommand{\hide}[1]{}

\begin{document}

\title{On the Asymptotic Behaviour of some Positive Semigroups }

\author{B. M. Makarow and M. R. Weber  }

\vspace*{1cm}
\institut{Sankt Petersburg State University and Technische Universit\"at Dresden  }

\preprintnumber{MATH-AN-09-2000 }
 
 \date{\small{}}

\sloppy
\makepreprinttitlepage 

\maketitle 

\setcounter{page}{1}


 \maketitle

\begin{abstract}
\noindent 
Similar to the theory of finite Markov chains it is shown that in a Banach space X ordered by a closed cone K with nonempty interior $\Int(K)$ 
a power bounded positive operator $A$ with compact power such that its trajectories for positive vectors eventually flow into $\Int(K)$, 
defines a "limit distribution", i.e. its adjoint operator has a unique fixed point in the dual cone. 
Moreover, the sequence $\{A^n\}_{n\in \mathbb{N}}$ converges with respect to the strong operator topology and for each functional $f\in X'$ 
the sequence $\{(A^*)^n(f)\}_{n\in \mathbb{N}}$ converges with respect to the weak*-topology (Theorem 5). 
If a positive bounded $C_0$-semigroup of linear continuous operators $\{S_t\}_{t\geq 0}$ on a Banach space contains a compact operator and the 
trajectories of the non-zero vectors $x\in K$ have the property from above then, in particular, $\{S_t\}_{t\geq 0}$ and $\{S^*_t\}_{t\geq 0}$ 
converge to their limit operator with repsect to the operator norm, respectively (Theorem 4). 
For weakly compact Markov operators in the space of real continuous functions on a compact topological space a corresponding result can 
be derived, that characterizes the long-term behaviour of regular Markov chains.
\end{abstract}

\section{Introduction}
%
%
The main purpose of our paper is to show that the method which is used to prove the
existence of a limit distribution in the theory of stationary Markov chains (see for
example \cite{Fel2}, chapt.VII,\S7, \cite{KemSn}, chapt.IV) can be transfered to a much more general situation.
The operator corresponding to the Markov chain is replaced by a positive
semigroup of operators acting in a Banach space ordered by a cone with nonempty interior, and the
condition of regularity of the Markov chain (in the sense of \cite{KemSn}) is 
transformed into the condition of strongly positivity of the operator or 
into an even more general condition (see condition 1) in the 
Theorems \ref{t1} --- \ref{t3}). 

In particular, the results generalize Theorem 1 of \cite{Fel1} and show that the limit distribution exists provided the
operator of the random walk on a compact space is weakly compact and satisfies the
condition of regularity. 

The main result (Theorem \ref{t2}) refers to the case, where a positive 
$C_0$-semigroup $\{S_t(x)\}_{t\geq 0}$ of linear continuous operators acts
in an ordered Banach space $X$. 
The semigroup of
operators is supposed to be uniformly bounded and to contain some compact operator. 
 Moreover, for each nonzero vector $x\in K$ its trajectory
$\{S_t(x)\}_{t\geq 0}$ eventually flows into $\Int(K)$. 
Then the following alternative takes place: either the operators $S_t$ for $t\to+\infty$ 
converge to $0$ with respect to the strong operator topology or the semigroup possesses a
common fixed point $u$ in $\Int(K)$, the adjoint operators 
$\{S^*_t\}_{t\geq 0}$ have a unique 
fixed point $f_0$ ("limit distribution") in the dual cone and finally,
$\{S_t\}_{t\geq 0}$ and $\{S^*_t\}_{t\geq 0}$ both for $t\to+\infty$ converge to the operator
$A_0=f_0\times u$ with respect to the norm topology.

Theorem \ref{t1} shows that the requirement of compactness can be 
considerably relaxed. However, in this case the convergence of the operators 
$\{S_t\}_{t\geq 0}$ and $\{S_t\}_{t\geq 0}$  to
$A_0$ takes place only with respect to the strong operator 
topology.

The Theorems \ref{t3} - \ref{t5} deal with sequences $\{A^n\}_{n\in \mathbb{N}}$ of iterates 
of a positive operator $A$ and are obtained as particular cases of the general 
results.

%
\section{Preliminaries}
%
We remember that a nonempty subset $K$ in vector space 
$X$ is a {\it wedge}, if $x, y\in K, \lambda,\mu\geq 0$ implies $\lambda x+\mu y\in K$. If in
addition $x, -x\in K$ implies $x=0$, then $K$ is a {\it cone}.   
In what follows we consider real normed spaces $(X,K,\|\cdot\|)$ in
which an order is introduced by means of a closed cone $K$.

\begin{defi}\label{D2}
{\rm Any complete normed space which is ordered by a closed cone  
we will call  an {\it ordered Banach space}
and denote it by $(X, K, \|\cdot\|)$.}
\end{defi}
Briefly we will write $X$ instead of $(X,K,\|\cdot\|)$ and denote its dual by $X'$. 
The closed ball in $X$ with radius $r>0$ and centered at the point $x$ is denoted by 
$B(x;r)$. We use the notations $x\in K$ and $x\geq 0$ synonymously.

A cone $K$ is said to be {\it generating} (or {\it reproducing}), if each vector
$x\in X$ has a representation as $x=x_1-x_2$, where $x_1,x_2\in K$.  

A cone $K$ is said to be {\it normal}, 
if there exists a positive number $\delta$, such that 
$\|x+y\|\ge\delta\max\{\|x\|,\|y\|\}$ for any $x,\ y \in K$.

A cone $K$ is said to be {\it nonflat},  
if there exists a positive constant $\gamma>0$, such that 
each element $x\in X$ is representable as $x=x_1-x_2$, 
where $x_i\in K$ and $\|x_i\|\leq \gamma \|x\|\; (i=1,2)$.

A linear functional defined on $X$ is said to be {\it positive}, if 
it takes on nonnegative values on all vectors of the cone $K$.
The set of all positive functionals of $X'$ is called the  
{\it dual wedge} and will be denoted by $K'$, i.e. 
$K'=\{f\in X\colon f(x)\geq 0 \ \mbox{for all} \  x\in K\}$. 

The following result goes back to M.G. Krein and V.L. \v{S}mulian (s. \cite{KLS})
\begin{theo}\label{t0}
If the cone $K$ is closed and normal then the wedge $K'$ is a closed 
genera\-ting cone, i.e. each 
functional of $X'$ has a representation as a difference of two 
positive functionals. Moreover, $K'$ is nonflat .
\end{theo}
\prf For the first part of the theorem see \cite{KLS}. We restrict ourselves to the proof 
of the nonflatness of $K'$. 
Let be $B_+^*=\{f\in K'\colon \|f\|\leq 1\}$ and $E=B_+^*-B_+^*$. According to the 
Banach - Alaoglu  Theorem (s. \cite{KA3},\ chapt.III, \S3) the set $B_+^*$ is
weak$^*$-compact, and therefore $E$ is closed. From the first part it follows that 
$X'=\bigcup_{n\in\mathbb{N}}nE$, and so $0$ is an interior point of $E$, i.e. for some
$r>0$ the ball $B^*(0;r)$ (in $X'$) belongs to $E$. This means
$r\frac{f}{\|f\|}\in E$ for each $f\in X',\ f\neq 0$, and implies that each 
functional $r\frac{f}{\|f\|}$ can be represented as $f=f_1-f_2$, where 
$f_1, f_2\in B_+^*$. Now $\gamma^*=\frac{1}{r}$ can be taken as the constant 
of nonflatness of the cone $K'$. \qed

For a convenient refering we list some more properties (s. \cite{Kra}, \cite{KLS}) 
of the space $X$, its dual $X'$, of  the cone $K$ and its dual cone $K'$ which are 
frequently used further on. 

\medskip
In the sequel we assume that the cone $K$
is closed and normal and satisfies $\Int(K)\neq\varnothing$. 
\begin{itemize}
\item[a)] {\it The cone $K$ is nonflat.} 

Indeed. Fix $u\in \Int(K)$. Then $u$ belongs to $K$ together with some
closed ball centered at $u$, i.e. $\overline{B}(u;r)\subset K$ for some $r>0$.
Then for any $x\in X$ 
\begin{equation}\label{f1}
\frac{\|x\|}{r}u \pm x \in K.
\end{equation}  
Put now $x_1=\frac{1}{2}\left(\frac{\|x\|}{r}u+x\right)$ and 
$x_2=\frac{1}{2}\left(\frac{\|x\|}{r}u-x\right)$. One has  $x_1, x_2\in 
K, \; x=x_1-x_2$ and $\|x_i\|\leq
\frac{1}{2}\left(\frac{\|u\|}{r}+1\right)\|x\|$. As the constant $\gamma$ (of 
nonflatness of $K$) can be taken the number $\frac{1}{r}\|u\|$. 
\item[b)] {\it Each linear positive functional $f$ on $X$ is continuous and satisfies 
the condition $f(x)>0$, if $f\in K',\ f\neq 0,\ x\in \Int(K)$.}

The relation (\ref{f1}) implies $\mp f(x)\leq
\frac{\|x\|}{r}f(u)$, which shows that $f$ is bounded on the unit
ball of $X$ and $\|f\|\leq \frac{1}{r}f(u)$. If $f\neq 0$ then $f(u)>0$.  
\item[c)] {\it  Each additive and positive homogeneous functional $f$ on $K$ with values in 
the nonnegative reals extends uniquely to a linear positive functional on the whole $X$.}
 
Indeed, if  
$x\in X$ is an arbitrary vector then $x=x_1-x_2$, where $x_i\in K$. 
Put 
\[  f(x)=f(x_1)-f(x_2).   \] 
It is easy to see that the functional $f$ is the required extension. 
We omit the standard proof (based on the nonflatness of $K$) of both the correctness 
of the definition and the uniqueness of the extension. 

\item[d)] {\it For any $x\in K, x\neq 0$ there exists a functional 
$f\in K'$ such that $f(x)>0$.}

Indeed, according to the theorem on a sufficient number of functionals there is a functional
 $f\in X'$ such that $f(x)\neq 0$. Since $f=f_1-f_2$ with $f_1, f_2  \in K'$, at least one 
 of the nonnegative  numbers $f_1(x), f_2(x)$ is strongly positive. 

\medskip
Remember that a set $D\subset K$ is called a {\it base of the cone} $K$, if $D$ is
convex and each vector $x\in K, \ x\neq 0$ has a unique representation as $x=\lambda y$, where 
$\lambda>0$ and $y\in D$. 

In the sequel we are interested in bases of the dual cone $K'$.
The existence of interior points in the cone $K$ guarantees that the cone 
$K'$ possesses a base. 
  
\item[e)] Let now $\mathcal{F}$ be an arbitrary base of the cone$K'$. Then the closedness 
of $K$ implies the following important property: If $x,y\in X $ then 
\[ 
 x\leq y  \quad \mbox{is equivalent to}\quad f(x)\leq f(y)\quad 
 \mbox{for all}\quad f\in \mathcal{F} \ \; .
\]
and consequently $x=y$ is equivalent to $f(x)=f(y)$ for all 
$f\in \mathcal{F}$. Moreover, together with b) one has $x\in \Int(K)$ if and only if
$f(x)>0$ for each $f\in \mathcal{F}$.
\item[f)]  For an arbitrary fixed element $u\in \Int(K)$ denote 
\[   \mathcal{F}=\mathcal{F}_u=\{f\in K'\colon f(u)=1\}.    \]
{\it Then the set $\mathcal{F}$ is bounded, weak*-compact and is a
base of the dual cone $K'$}. 

The relation (\ref{f1}) implies the estimate 
\begin{equation}\label{f02}  
|f(x)|\leq \frac{1}{r}\|x\|  \quad \mbox{for any} \quad f\in \mathcal{F},\ x\in
X.
\end{equation}
Because of its weak$^*$-closedness the set $\mathcal{F}$ is weak$^*$-compact by the
Banach - Alaoglu  Theorem.  
The set $\mathcal{F}$ is convex, and property b) implies that $f(u)>0$ for $f\in K,\ f\neq
0$. Therefore, $\mathcal{F}$ is a base of the dual cone.

\medskip
 By means of the interior point $u$ of the cone $K$ one can define
 the following nonnegative functional on $X$
\[
 \|x\|_u = \inf\{\lambda\geq 0\colon -\lambda u\leq x\leq \lambda
u\}
\]
which is called the {\it u-norm}.  
Notice that the $u$-norm of an element $x$ can be calculated also by the 
formula 
\[ \|x\|_u =\sup\{|f(x)|\colon f\in \mathcal{F}\}. 
\] 
It is clear that the u-norm is actual a norm and that it is monotone on $K$, i.e. $x\leq y$ implies
$\|x\|_u\leq \|y\|_u$. 
\item[g)] {\it The u-norm is equivalent to the original norm on $X$.}
Indeed, (\ref{f02}) implies  
\[ 
\|x\|_u=\sup\{|f(x)|\colon f\in \mathcal{F}\}\leq \frac{1}{r}\|x\|.
\] 
On the other hand, let $x\in X, \ f\in X',\ \|f\|=1$ and $f(x)=\|x\|$. Then
$f=f_1-f_2$, where $f_i\in K'$ and $\|f_i\|\leq \gamma^*\|f\|$ and 
$\gamma^*$ denotes the constant of nonflatness of the cone $K'$. Then 
\[ \|x\|=f(x)\leq |f_1(x)|+|f_2(x)|
	 \leq f_1(u)\|x\|_u+f_2(u)\|x\|_u\leq  
 C_0\|x\|_u,  
\]
where $C_0=2\gamma^*\|u\| $. 

\end{itemize}
Summing up we have that for each vector $u\in \Int(K)$ there is a constant $C_u>0$ such
that for each $x\in X$ 
\begin{equation}\label{f05} 
C_u^{-1} \|x\|\leq \|x\|_u\leq C_u\|x\|.
\end{equation}

\bigskip 
We consider now  positive operators on $(X, K,\|\cdot\|)$.  
By $L(X)$ we denote the vector space of all
linear continuous operators on $X$, equipped with the usual norm and the
order. For $A\in L(X)$ we write $A \geq 0$ if $A(K)\subset K$. Such operators we 
will call {\it positive}.  The simple
properties of such operators are gathered in the  
\begin{lem}\label{l1}
Let $(X,K,\|\cdot\|)$ be an ordered normed real vector space,
and let $A$ be a positive linear continuous operator on $X$. 
Assume there exists a vector $u\in \Int(K)$ such that $A(u)=u$. 
Let $\mathcal{F}_u=\{f\in K'\colon f(u)=1\}$ be the base of $K'$
corresponding to the vector $u$. 

Then the following statements hold. 
\begin{itemize}
\item[(i)] The adjoint operator $A^*$ is positive.
\item[(ii)] $A^*(\mathcal{F}_u)\subset\mathcal{F}_u$. 
\item[(iii)] $\|A(x)\|_u\leq \|x\|_u$ for each $x\in X$.
\item[(iv)] $\{A^n\}_{n\in \mathbb{N}}$ is a norm bounded sequence in $L(X)$.
\item[(v)] $A(x)\in \Int(K)$ for each $x\in\Int(K)$. 
\end{itemize}
\end{lem}
\prf 
(i) For an arbitrary vector $f\in K$ we have
$A^*(f)(x)=f(A(x))\geq 0$ for any $x\in K$. Thus $A^*(f)$ belongs
to $K'$.

\medskip
(ii) If $f\in \mathcal{F}_u$ then by (i) $A^*(f)\in K'$. Since $A(u)=u$ and
$A^*(f)(u)=f(A(u))=f(u)=1$ we obtain $A^*(f)\in \mathcal{F}_u$. 

\medskip
(iii) From (ii) follows that for each $x\in X$ one has 
\[  \|A(x)\|_u=\sup_{f\in \mathcal{F}_u}|f(A(x))|= \sup_{f\in \mathcal{F}_u}|(A^*(f))(x)|
\leq \sup_{f\in \mathcal{F}_u}|f(x)|= \|x\|_u .
\]

\medskip
(iv) It is convenient to use the $u$-norm in $X$ which, as was
shown in property g), is equivalent to the norm $\|\cdot\|$. Use now the inequality
(\ref{f05}) and (iii) then  
\[  
\|A^n(x)\| \leq  C_u\|A^n(x)\|_u \leq C_u\|A^{n-1}(x)\|_u 
 \leq  \ldots  \leq C_u\|x\|_u \leq C_u^2 \|x\| \quad  \mbox{for all}\quad  x\in X.
\]
Therefore, $\|A^n\|\leq C_u^2$ for all $n\in \mathbb{N}$. 

\medskip
(v) If $x\in \Int(K)$ then 
in view of (ii) one has $f(A(x))=(A^*(f))(x)>0$ for each $f\in \mathcal{F}_u$. 
By property e) it follows $A(x)\in \Int(K)$. 
\qed

We need also the following auxiliary result concerning positive operators
\begin{lem}\label{l2}
Let $A$ be a positive operator which satisfies the following
conditions  
\begin{enumerate}
\item[1)] $A(\Int(K))\subset \Int(K)$;
\item[2)] for each vector $x\in K, \ x\neq 0$  there exists a natural  
$n_x$ such that $A^{n_x}(x)\in \Int(K)$. 
\end{enumerate}  
Then for any compact set $R\subset K$ such that $0\notin R$ 
there is a natural number $p$ with $A^p(R)\subset \Int(K)$. 
\end{lem}
\prf  The condition  1) implies 
\begin{equation}\label{f21}
(A^n)^{-1}(\Int(K))\subset(A^{n+j})^{-1}(\Int(K)) 
\end{equation} 
for all $n,j\in\mathbb N$.

In view of condition 2) for each vector $z\in R$ there exists a power 
$n_z$ such that $A^{n_z}(z)\in \Int(K)$. Therefore the sets 
$(A^{n_z})^{-1}(\Int(K))$ form an open covering of $R$.  
Consider any finite subcovering 
\[
(A^{n_{z_1}})^{-1}(\Int(K)), \ (A^{n_{z_2}})^{-1}(\Int(K)),
     \ldots,(A^{n_{z_s}})^{-1}(\Int(K))
\]
and let be $p=\max\{n_{z_1}, n_{z_2},\ldots,n_{z_s}\}$.
Then taking into consideration inclusion
(\ref{f21}) the family consisting of $s$ exemplars of $(A^p)^{-1}(\Int(K))$
also covers the set $R$. Actually we have 
$R\subset (A^p)^{-1}(\Int(K))$.
This shows that $A^p(z)\in \Int(K)$ for each $z\in R$.
\qed

\medskip
\textbf{Remark 1} The condition 1) of Lemma \ref{l2} is fulfiled, if the operator $A$ is
positive and possesses a fixed point $u$ such that $u\in \Int(K)$ (s.
Lemma \ref{l1}(v)).
%

%

\section{Main results}
%
We need the following notations (\cite{CleHei}). 
\begin{defi}\label{D1}
{\rm Let $X$ be a (real) Banach space. A family $\{S_t\}_{t\geq 0}$ of operators in
$L(X)$ is called a} one-parameter semigroup of bounded linear operators 
{\rm  if
$S_0=I, \; S_{s+t}=S_s S_t \; (s,t\geq 0)$, where $I$ denotes the identity
operator on $X$. 

\smallskip
If, in addition, the function $t\mapsto S_t$ is continuous with respect to the
strong operator topology, i.e. the function $t\mapsto S_t(x)$ is norm-continuous on
$[0, +\infty)$ for each $x\in X$, then $\{S_t\}_{t\geq 0}$ is called a} {\it
strongly continuous semigroup}, {\rm or also a} $C_0$-semigroup. 

\smallskip
{\rm A $C_0$-semigroup $\{S_t\}_{t\geq 0}$ in an ordered Banach space is called} positive, 
{\rm if each operator $S_t$ is positive ($t\geq 0$).}
\end{defi}

\medskip
The results we are going to prove are valid in real ordered Banach
spaces $(X,K,\|\cdot\|)$, briefly denoted by $X$, where $K\subset X$ 
is a closed normal cone which satisfies the condition $\Int(K)\neq \varnothing$.  

\begin{theo}\label{t1}
Let $(X,K,\|\cdot\|)$ be an ordered Banach space and $\{S_t\}_{t\geq 0}$ a 
positive $C_0$-semigroup of 
operators in $L(X)$ which satisfies the following conditions  
\begin{itemize}
\item[1)] for each vector $x\in X,\,x\neq 0$ there exists a
 number $t_x\in [0,+\infty)$ such that $S_{t_x}(x)\in \Int(K)$; 
\item[2)] for each vector $x\in K$ its trajectory
$\{S_t(x)\}_{t\geq 0}$ is relatively compact; 
\end{itemize}
Then the family $\{S_t\}_{t\geq 0}$ converges pointwise for $t \to +\infty$
to some operator $A_0$. 

If that operator $A_0$ is not the zero one, then there exist a vector 
$u\in \Int(K)$ and a functional $f_0\in K', f_0\neq 0$ such that 
\begin{itemize}
\item[(i)] $S_t(u)=u$, $S^*_t(f_0)=f_0$ for any $t\geq 0$ and moreover $f_0(x)>0$ if 
$x\in K,\,x\neq 0$; 
\item[(ii)] $A_0=f_0\otimes u$ 
\item[(iii)] for each $f\in X'$ one has 
$S^*_t(f)\konv_{t\to\infty} A_0^*(f)$ 
 with respect to the weak*-topology $\sigma(X',X)$;
\item[(iv)] $\lambda=1$ is a simple eigenvalue of the operators
$S_t$ and $S^*_t$ for all $t>0$.
\end{itemize}
\end{theo}

\prf I. First of all we show that all operators $S_t$ for $t\geq 0$ have a common 
fixed point in $\Int(K)$, provided the family $\{S_t\}_{t\geq 0}$ does not pointwise
converge to the zero operator ${\mathbf 0}$ for $t\to +\infty$. 
According to the principle of uniform boundedness the condition 2) implies that the norms of all operators of the semigroup 
$\{S_t\}_{t\geq 0}$ are bounded, i.e. there exists a constant $C$ such
that $\|S_t\|\leq C$ for all $t\in [0,\infty)$. 
If $\{S_t\}_{t\geq 0}$ does not converge to the zero operator, then 
\begin{equation}\label{fh2}
\limsup\limits_{t\to+\infty} \, \|S_t(x_0)\|>0.
\end{equation}
holds for some vector $x_0\in X$. Since $K$ is generating we may assume
$x_0\in K$. Now we show that $\inf\limits_{t\geq 0}\|S_t(x_0)\|>0$.
Indeed, if the contrary is assumed we find an increasing sequence
$\{s_k\}_{k=1}^\infty \subset [0,+\infty), \; s_k\to +\infty$ such that 
$\|S_{s_k}(x_0)\|\konv_{k\to\infty} 0$. Then an arbitrary sequence
$\{t_n\}_{n=1}^\infty$ with $t_n\to +\infty$ satisfies $s_{k_n}\leq t_n\leq
s_{k_{n+1}}$, where $s_{k_n}\konv_{n\to\infty}+\infty$ and therefore, 
\[   
     \|S_{t_n}(x_0)\|=\|S_{t_n-s_{k_n}}S_{s_{k_n}}(x_0)\| \leq
     C\|S_{s_{k_n}}(x_0)\| \konv_{n\to\infty} 0 .
\]
This contradicts to the inequality (\ref{fh2}).

Consequently, the points of the trajectory $\{S_t(x_0)\}_{t\geq 0}$ belong
to $K$ and their norms are separated from zero. 
Denote by $Q_0$ the closure of that trajectory. Then $Q_0$ is compact and
$S_t(Q_0)\subset Q_0$ for any $t\geq 0$. Denote the closure of the
convex hull of the set $Q_0$ by $Q_1$.
We show now that the zero-vector also does not belong to $Q_1$. 
From $Q_0\subset K$ and $0\notin Q_0$ we find (property d))
\[
 Q_0\subset\bigcup_{f\in K'}\{x\in X \colon f(x)>0\}.
\]
By selecting a finite covering we have 
\[
 Q_0\subset\bigcup_{k=1}^N\{x\in X : f_k(x)>0\},
\]
where $f_k\in K'$. Put $g=\sum_{k=1}^Nf_k$. Then, obviously, 
$Q_0\subset \{x\in X : g(x)>0\}$ and therefore, there is some $\sigma>0$ 
such that the set $Q_0$ is contained in the closed convex set 
$K\cap\{x\in X \colon g(x)\ge\sigma\}$, 
which in turn does not contain zero. Now it is clear that $0$ does not
belong to the closed convex hull $Q_1$ either. Since $S_t(Q_1)\subset Q_1$ 
 for any $t\ge0$ and since the operators of the family
$\{S_t\}_{t\geq 0}$ commute in pairs we are able to apply the
Markov-Kakutani Theorem to that family of operators on $Q_1$ 
(s.\cite{Edw}, chapt.III.3.2) and to conclude that they possess  
a common fixed point, say $u$, in  $Q_1$. It is obvious  that $u\neq 0$ and in
view of condition 1) $u\in \Int(K)$. This completes the first step of the proof
and allows us to apply the Lemmata \ref{l1} and \ref{l2} to each of the 
operators $S_t$. In particular, from Lemma \ref{l1}(v) we get 
 $S_t(\Int(K))\subset \Int(K)$ for any $t>0$, and from condition 1) there follows
 that 
 \[ S_s(x)=S_{t_x+t}(x)=S_t(S_{t_x}(x))\in \Int (K)\quad \mbox{for any } s\geq
 t_x ,\]
 i.e. each trajectory starting at $x\in K$ will stay eventually in $\Int(K)$. 
 
II. We prove now that  for each vector $x\in K$  there exists the limit 
$\lim\limits_{t\to\infty}S_t(x)$ with respect to the norm.

We introduce the sets 
\[ 
              \mathcal{F}=\{ f\in K'\colon f(u)=1 \}
\] 
and  $\mathcal{S}_+=\{ x\in K \colon \max_{f\in\mathcal{F}}f(x)=1 \}$. As was mentioned
above (property f)) $\mathcal{F}$ is a convex weak$^*$-compact set which 
is a base of the cone $K'$. 
Both sets  $\mathcal F$ and $\mathcal{S}_+$ are closed with respect to the norm, and 
$u\in \mathcal{S}_+, \ 0\notin \mathcal{S}_+$. Observe that $S_t(u)=u$ for any $t>0$ implies  
$S^*_t(\mathcal F)\subset\mathcal F$ (Lemma \ref{l1}(ii)).

For any vector $x \in X$ and $t\in[0,+\infty)$ we define the numbers  
\[
 M_x(t)=\sup_{f\in \mathcal{F}}f(S_t(x)) \quad \mbox{and} \quad  
 m_x(t)=\inf_{f\in \mathcal{F}}f(S_t(x)) .
\] 
The equation $f(S_{s+t}(x))=\left(S_s^*(f)\right)(S_t(x))$ and the
inclusion $S_s^*(\mathcal{F})\subset \mathcal{F}$  imply 
\begin{equation}\label{f22}
m_x(t)\leq m_x(t+s)\leq M_x(t+s)\leq M_x(t)\quad \mbox{for all}\quad s,t\geq 0. 
\end{equation} 

Therefore the functions $M_x(t)$ and  $m_x(t)$ are monotone and possess 
finite limits at infinity. 

Moreover, since for any $f\in \mathcal{F}$ the
inequalities  $m_x(t)\leq f(S_t(x))\leq M_x(t)$ can be written as
\[  f(m_x(t)\ u)\leq f(S_t(x))\leq f(M_x(t)\ u), \]
the inequality
\begin{equation}\label{fh3}
     m_x(t)\ u \leq S_t(x)\leq M_x(t)\ u
\end{equation}
holds for each $x\in K$ (see property e)).

The main aspect of the proof is to establish the relation 
\begin{equation}\label{f22b}
\delta_x(t)\equiv M_x(t)-m_x(t)\konv_{t\to+\infty} 0
\end{equation}
for each $x\in X$.

In view of (\ref{f22}) it suffices to prove that some sequence
$\{\delta_x(kt_0)\}_{k\in \mathbb{N}}$ for $t_0>0$ converges to $0$. Assume by way of 
contradiction that there is some  $x_0\in X$ such that 
$\delta_{x_0}(t)\not\rightarrow 0$ for $t\to+\infty$. Due to the
monotony of  $\delta_{x_0}(t)$ this means
\begin{equation}\label{fh4}
\delta_{x_0}(t) = M_{x_0}(t) - m_{x_0}(t)\ge\varepsilon
\end{equation}
for some $\varepsilon>0$ and all $t\in [0,\infty)$.

Consider now the set
\[
R_0=\left\{\frac{S_t(x_0)-m_{x_0}(t)u}{\delta_{x_0}(t)}, \
           \frac{M_{x_0}(t)u-S_t(x_0)}{\delta_{x_0}(t)} \colon \quad
t\in [0,\infty) \right\}.
\]
In view of condition 2), the inequality (\ref{fh4}) and the
boundedness of the functions $M_x(t)$ and  $m_x(t)$, the set $R_0$ turns out to be
relatively compact. Moreover, it is easy to see that $R_0\subset \mathcal{S}_+$.
Therefore the closure $R$ of $R_0$ is also contained in $\mathcal{S}_+$ and, in particular,
$0$ does not belong to $R$.

According to the Lemma \ref{l2} and the Remark 1 there is a natural number $p$
such that the compact set $Q=A^p(R)$ belongs to $\Int(K)$, where
$A=S_1$.

The bilinear form $\langle z,f\rangle=f(z)$ is strongly positive on the compact set
$Q\times \mathcal{F}$, where $Q$ is considered with the norm topology (induced from $X$) 
 and $\mathcal{F}$ with the weak*-topology (s. property f)). 
 The inequalities 
 \begin{eqnarray*}
|\langle z,f\rangle -\langle y,g\rangle| &\leq &|\langle z-y,f\rangle| +
|\langle y,f-g\rangle| \\
 & \leq &  \|z-y\|\ \|f\| + |\langle y,f-g\rangle|
\end{eqnarray*}
show that the bilinear form is continuous on the set $Q\times \mathcal{F}$. 
Therefore there is some positive number $\beta$ such that $f(z)>\beta$ for all 
$z\in Q$ and $f\in\mathcal{F}$. 
We shall assume $\beta<\frac{1}{2}$.

The vectors (remember that $A=S_1$)
\[
 A^p\left(\frac{A^n(x_0)-m_{x_0}(n)u}{M_{x_0}(n)-m_{x_0}(n)}\right)\quad
\mbox{and} \quad
A^p\left(\frac{M_{x_0}(n)u-A^n(x_0)}{M_{x_0}(n)-m_{x_0}(n)}\right),
\]
belong to $Q$ and, consequently, for each $f\in \mathcal F$ we have 
\[
 f\left(A^p\left(\frac{A^n(x_0)-m_{x_0}(n)u}
{M_{x_0}(n)-m_{x_0}(n)}\right)\right)\ge\beta
\quad  \mbox{and} \quad
f\left(A^p\left(\frac{M_{x_0}(n)u-A^n(x_0)}{M_{x_0}(n)-m_{x_0}(n)}\right)\right)
\ge \beta.
\]
This together with $\beta=f(\beta u)$ implies by e) 
\begin{eqnarray*}
A^{n+p}(x_0) & \geq & m_{x_0}(n)\,u+\beta(M_{x_0}(n)-m_{x_0}(n))\,u
\quad \text{ and } \\ 
A^{n+p}(x_0) & \leq & M_{x_0}(n)\,u-\beta(M_{x_0}(n)-m_{x_0}(n))\,u \quad
\mbox{for any}\quad n\in\mathbb N.
\end{eqnarray*}
Put now $n=kp$. Then
\[
 M_{x_0}((k+1)p)-m_{x_0}((k+1)p)\le
(1-2\beta)(M_{x_0}(kp)-m_{x_0}(kp)),
\]
and therefore
\[
 M_{x_0}(kp)-m_{x_0}(kp)\le (1-2\beta)^k(M_{x_0}(0)-m_{x_0}(0))\konv_{k\to\infty} 0.
\]
However this contradicts to (\ref{fh4}). So the relation (\ref{f22b}), i.e.
$M_x(t)-m_x(t)\konv_{t\to\infty} 0$, is proved.

\medskip
III. In order to complete the final part of the proof 
we denote for each $x\in X$
\[  f_0(x)= \lim_{t\to\infty}m_x(t).\]
From the inequalities 
\begin{equation}\label{f22d}
m_x(t)\, u \leq f_0(x)\, u \leq M_x(t)\, u 
\end{equation}
and (\ref{fh3}) by means of passing to the limit we obtain for each $f\in
\mathcal F$ 
\[ 
      f_0(x)=\lim_{t\to\infty}f(S_t(x))
\]
and so $f_0$ is an additive, homogeneous and nonnegative functional on
$X$ such that $f_0(u)=1$. 
The inequalities (\ref{fh3}) and (\ref{f22d})  further imply for $x\in K$
\begin{equation}\label{fh8} -\left(M_x(t)-m_x(t))\right)\,u\le S_t(x)-f_0(x)\,u\le
\left(M_x(t)-m_x(t))\right)\,u.
\end{equation}
Define now the rank-one operator $A_0$ by $A_0=f_0\otimes u$, i.e. 
$A_0(x)=f_0(x)\,u$  for $x\in X$.
From (\ref{fh8}) and (\ref{f05}) it follows that 
\[  
\|S_t(x)-A_0(x) \| \le C_u \|S_t(x)-A_0(x) \|_u \le
C_u\left(M_x(t)-m_x(t)\right)\konv_{t\to\infty} 0 . 
\]
This proves the statement (ii) of the theorem.

We finalize the proof of the statements (i) and (iii). Since
$S_{(n+1)t}(x)=S_{nt}(S_t(x))$ for each $x\in X, \ t>0$ and
$n\in \mathbb N$, after passing to
the limits as $n\to\infty$ we obtain $f_0(x)\,u=f_0(S_t(x))\,u$
which shows that
$f_0(x)=\big(S^*_t(f_0)\big)(x)$ for each $x\in X$, i.e.
$f_0=S^*_t(f_0)$ for $t>0$.
We show that $f_0(x)>0$ if $x\in K, x\neq 0$. For such $x$ there is some $t_x$
with $S_{t_x}(x)\in \Int(K)$. In view of property b) and the weak$^*$-compactness 
of $\mathcal{F}$ we get  $m_{x}(t_x)>0$. Then 
$f_0(x)=\lim\limits_{t\to\infty}m_x(t)\geq m_x(t_x)>0$.

\smallskip
From the already proved statement (ii) it follows that for each functional
$f\in\mathcal{F}$ the family $\{S^*_t(f)\}$ converges to $f_0$ with respect
to the weak$^*$-topology, i.e. 
\[  S_t^*(f)(x)=f(S_t(x))\konv_{t\to\infty} f_0(x) \quad \mbox{for each}\quad x\in X.  \]
It remains to notice that due to the facts that any functional $f\in X'$ is
representable as a difference of two nonnegative functionals (s. Theorem \ref{t0})
and that $\mathcal{F}$ is a base of the dual cone $K'$, the last relation
implies 
\[  
S_t^*(f)\konv_{t\to\infty} f(u)f_0\quad\mbox{for each}\quad f\in X'  
\]
with respect to the weak$^*$-topology. Now (iii) is proved.

\medskip
It remains to prove (iv), i.e. that $\lambda=1$ is a simple eigenvalue of the
operators $S_t$ and $S^*_t$ for $t>0$. Indeed, if $u'$ is another fixed point of 
$S_t$, then 
$S_{nt}(u')=u'$ for any $n\in\mathbb{N}$. Since
$S_{nt}(u')\konv_{n\to\infty}f_0(u') u$ 
one immediately has $u'=f_0(u')u$. 
That means the eigenspace of the operator $S_t$, corresponding to the eigenvalue $\lambda=1$, 
is one-dimensional. A similar argument shows the statement for the adjoint operator.   
\qed 

\begin{corr}\label{C1} 
For the operators $S_t$ and $A_0$ for each $t\in [0,\infty)$
and $n\in\mathbb N$  there hold the following relations
\begin{itemize}
\item[a)] $A_0^n=A_0$; 
\item[b)] $S_tA_0=A_0S_t=A_0$; 
\item[c)] $(S_t-A_0)^n=S_{nt}-A_0$. 
\end{itemize} 
\end{corr}
\prf a) - c) are obtained by a simple calculation which we will omit. \qed

\medskip
If the condition 2) of the theorem is replaced by a slighty stronger one, then the operators 
$S_t^*$, for $t\to \infty$, converge to the operator $A_0^*$ not only in the weak operator toplogy but also pointwise. The new
condition is well known in the theory of Markov chains (see, for example, \cite{nev}(Lemma
V.3.1).  
We come now to one of our main results.
\begin{theo}\label{t1a}
Let  $(X,K,\|\cdot\|)$ be an ordered Banach space and $\{S_t\}_{t\geq 0}$ a 
positive $C_0$-semigroup of operators in $L(X)$ which satisfies the following 
conditions   
\begin{itemize}
\item[1)] for each vector $x\in K,\,x\neq 0$ there exists a
 number $t_x\in [0,\infty)$ such that $S_{t_x}(x)\in \Int(K)$; 
\item[2)] there exist a number $\tau>0$ and a compact operator $V$ such that
$\|S_{\tau}-V\|<1$;
\item[3)] $\sup\limits_{t\ge 0}\|S_t\|<\infty$.
\end{itemize}
Then all statements  of Theorem \ref{t1} are valid. Moreover,  
\begin{itemize}
\item[(v)] the operators $S_t^*$, for $t\to\infty$, converge pointwise to the 
operator $A_0^*$.
\end{itemize}
\end{theo}
\prf 
First of all we prove that the condition 2) of Theorem \ref{t1} is satisfied. If $x\in X$
then it suffices to show that for an arbitrary fixed $\varepsilon >0$ the trajectory 
$\{S_t(x)\}_{t\ge 0}$
possesses a relatively compact $\eps$-net. Put $C=\sup\limits_{t\geq 0}\|S_t\|$, 
$W=S_{\tau}-V$ and $q=\|W\|$. Obviously $C<\infty$ and 
\[
S_{n\tau}=S_{\tau}^n=W^n+V_n,\ \] 
where $V_n$ is some compact operator and $\|W^n\|\leq q^n\konv_{n\to \infty} 0$.

Fix a sufficiently large $N$ such that $q^N<\eps$.  Notice that the trajectory 
$\{S_t(x)\}_{t\geq 0}$ is contained in the closed ball $B_x=B(0;C\|x\|)$. If $t=N\tau+t', t'>0$
then $S_{t'}(x)\in B_x$, and therefore 
\[
  \|S_t(x)-V_N(S_{t'}(x))\| =\|(S^N_{\tau}-V_N)(S_{t'}(x))\|\leq q^NC\, \|x\|<C\eps \|x\|.
\]
Now it is immediate that the relatively compact set 
\[  \{S_t(x)\colon 0\leq t\leq N\tau\}\cup V_N(B_x) \] 
is a $C\eps\|x\|$-net for the trajectory $\{S_t(x)\}_{t\geq 0}$, and so the condition 2) of
the Theorem \ref{t1} holds. 

We prove now the statement ($v$). Assume first $A_0\neq 0$. In this case it suffices to show
$\|S_t^*(f)-f_0\|\konv_{n\to\infty} 0$ for each $f\in\mathcal{F}$. 
Let $\eps$ be an arbitrary positive number and $N, t', V_N$ be the same as above. 
Let be $H=V_N(B(0;1))$. In view of the equalities  
\[  f(S_{N\tau}(x))=f(V_N(x))+f(W^N(x)), \quad
f_0(x)=(S_{N\tau}^*(f_0))(x)=f_0(V_N(x))+f_0(W^N(x)) \]  
we obtain for $t=N\tau+t'$ the estimate 
\begin{eqnarray*}
\|S_t^*(f)-f_0\|& \leq & \sup_{\|x\|\leq 1}|S_{t'}^*(f)(V_N(x))-f_0(V_N(x))| \\ 
                &  +   & \sup_{\|x\|\leq 1}|S_{t'}^*(f)(W^N(x))|+ \sup_{\|x\|\leq 1}|f_0(W^N(x))| \\
                & \leq & \sup_{y\in H}|S_{t'}^*(f)(y)-f_0(y)|+C\|f\|+\|f_0\|q^N \\
		&   <  & \sup_{y\in H}|S_{t'}^*(f)(y)-f_0(y)|+(C\|f\|+\|f_0\|)\eps. 
\end{eqnarray*} 
This estimate holds for any $t'>0$. 
Because of $S_{t'}(f)\konv_{t'\to\infty} f_0$ with respect to the weak$^*$ topology and the
relative compactness of the set $H$ the supremum at the right side of the inequality
converges to $0$ if $t'\to\infty$ by the the theorem on uniform convergence on compact
sets. Consequently, for sufficiently large $t'$ we obtain 
\[  \|S_t^*(f)-f_0\|<\eps +(C\|f\|+\|f_0\|)\eps, \]
what has to be shown. 

If $A_0=0$ then for each $f\in X'$ the given proof is applicable if $f_0=0$ is assumed. 
\qed

\begin{corr}\label{C2} 
Under the conditions of Theorem \ref{t1a} the operator 
$T:=T(t)=I-S_t+A_0$ is invertible for any $t>0$ and 
\begin{equation}\label{f30}   
T^{-1}=I + \sum_{n=1}^\infty(S_{nt}-A_0),  
\end{equation}
where the series converges pointwise.   
\end{corr}
\prf
We show first that $\ker(T)=\{0\}$ for all $t>0$.
(s. \cite{Edw}, propositions 9.10.2, 9.10.5).

Assume $A_0\neq {\bf 0}$. If $T(x_0)=0$ then apply 
the  operator $A_0$ to the equation $-A_0(x_0)=x_0-S_t(x_0)$ 
and by taking into consideration the statements a) and b) of Corollary \ref{C1} we see that 
$A_0(x_0)=f(x_0)u=0$. Therefore $f_0(x_0)=0$, and due to statement (i) of the theorem we get $x_0=0$.

If  $A_0={\bf 0}$, then the operator $S_t$ can not have any nonzero 
fixed point $x_0$, since in the opposite case there would be 
$S_{nt}(x_0)=x_0\not\to A_0(x_0)=0$. So, 
\[ 0=T(x_0)=x_0-S_t(x_0)+A_0(x_0)=x_0-S_t(x_0)   \] 
means $S_t(x_0)=x_0$ and implies that the kernel of the operator $T$ is trivial for any $t>0$. 

\medskip  
The proof of invertibility of the operator $T$ for all $t>0$ now follows. 
By keeping the notation of the theorem we put $W=S_{\tau}-V, \ q=\|W\|$. Put also 
$U\equiv S_t-A_0$. Notice that according to the Corollary \ref{C1} and Theorem \ref{t1}
the sequence $U^n=S_{nt}-A_0$ converges to $0$ pointwise. 

We fix now some $m\in\mathbb{N}$ such that  $2C\,q^m<1$, where 
$C=\sup_{t\ge 0}\|S_t\|$, and prove that $T$ is invertible for 
$t>m\tau$. Remember that $S_{m\tau}=W^m+V_m$, where $V_m$ is some compact 
operator. If $t=m\tau+\sigma, \ \sigma\ge 0,$  then $T$ (at the moment $t$) is equal to 
\[
T=I-S_{m\tau}(S_\sigma-A_0)=I-W^m(S_\sigma-A_0)-V_m(S_\sigma-A_0).
\]
Notice that the operator $R=I-W^m(S_\sigma-A_0)$ in invertible because of 
$\|W^m(S_\sigma-A_0)\|\le q^m(C+\|A_0\|)\le 2C\,q^m<1$. Hence 
\[
T=R\big(I-R^{-1}V_m(S_\sigma-A_0)\big).
\]
At the same time the operator $I-R^{-1}V_m(S_\sigma-A_0)$ is invertible since, 
in view of $\ker{T}=\{0\}$, its kernel is trivial, and the  operator 
$R^{-1}V_m(S_\sigma-A_0)$ is compact together with $V_m$. 
The invertibility in this case of the operator $T$ is established.

In the case of $0<t<m\tau$, we use the identity 
\begin{equation}\label{f31}
I-U^n=T(I+U+\dots +U^{n-1}).
\end{equation}
If $nt>m\tau$ then by what has been shown above the operator $I-U^n=I-S_{nt}+A_0=T(nt)$ is 
invertible. It follows from (\ref{f31}) that also the operator $T=T(t)$ is invertible 
(since due to the invertibility of the operator $I-U^n$ it shall be injective and 
surjective). 

From the identity (\ref{f31}) and the invertibility  of the  operator $T$ it follows 
that
\begin{equation}\label{f32}
T^{-1}-T^{-1}U^n=I+U+\dots +U^{n-1}.
\end{equation}
Since  $U^n\to 0$ pointwise this proves that the decomposition (\ref{f30}) takes place.
\qed

\medskip
Our next result is 
\begin{theo}\label{t2}
Let  $(X,K,\|\cdot\|)$ be an ordered Banach space and $\{S_t\}_{t\geq 0}$ a 
positive $C_0$-semigroup of 
operators in $L(X)$ which satisfies the following conditions   
\begin{itemize}
\item[1)] for each vector $x\in K,\,x\neq 0$ there exists a
 number $t_x\in [0,\infty)$ such that $S_{t_x}(x)\in \Int(K)$; 
\item[2)] for some $\tau>0$ the operator $S_\tau$ is compact;
\item[3)] $\sup\limits_{t\ge 0}\|S_t\|<\infty$.
\end{itemize}
Then the statements
of Theorem \ref{t1} are valid. Moreover,  
\begin{itemize}
\item[(v)] $\|S_t-A_0\|\konv_{t\to\infty}0$, i.e. the operators $S_t$ (and, of
course, the adjoint operators)  
converge to  $A_0$ (to $A_0^*$) with respect to the norm.
\end{itemize}
\end{theo}

\prf Since the condition 2) of this theorem is stronger than the corresponding condition
2) of the previous theorem and the other ones coincide, it is left to prove only
statement (v).
\hide{s with them condition 
of Theorem \ref{t1} we check that also the condition 2) of 
Theorem \ref{t1} is satisfied. 
Let be $\|S_t\|\le C$ for all $t\in[0,+\infty)$.
Due to $\|S_t(x)\|\le C\,\|x\|$ the vectors $S_t(x)=S_s S_{t-s}(x))$ belong to
the image under the compact operator $S_s$ of the closed ball $B(0;
C\,\|x\|)$ for $t>s$.
Therefore, taking into consideration the property of an $C_0$-semigroup, the trajectory of an arbitrary point $x$ belongs to the relatively
compact set 
$\{S_t(x)\colon 0\le t\le s\}\cup S_s(B(0; C\,\|x\|))$ 
what ensures that the condition 2) of Theorem \ref{t1} holds.
We prove now the statement (v). }

Due to b) of Corollary \ref{C1} for $t>\tau$ one has 
\[
       S_t-A_0=(S_{t-\tau}-A_0)\,S_\tau,
\]
and if $Q$ denotes the closure of the image under $S_\tau$ of the unit ball, 
therefore 
\[
\|S_t-A_0\|=\sup_{\|x\|\le 1}\|(S_{t-\tau}-A_0) S_\tau(x))\|=
\sup_{y\in Q}\|(S_{t-\tau}-A_0)(y)\|,
\]
Since $S_{t-\tau}\to A_0$ pointwise as $t\to \infty$ and $Q$ is compact one has (s.
\cite{Bour1},chapt.III \S3, prop.5)
\[
\sup_{y\in Q}\|(S_{t-\tau}-A_0)(y)\| \ \konv_{t\to+\infty} \ 0.
\]
This completes our proof. \qed 

\medskip

\begin{corr}\label{C3} 
Under the conditions of Theorem \ref{t2} the operator  
$T:=T(t)=I-S_t+A_0$ 
is invertible for $t>0$ and
\begin{equation}\label{f23}
    T^{-1}(t) = I+\sum_{n=1}^\infty (S_{nt}-A_0).
\end{equation}
where the series converges with respect to the norm.
\end{corr}
\prf
The invertibility of the operator $T$ has been proved in Corollary \ref{C2}. 
Therefore it remains to pass to the limit  in the identity (\ref{f32}) by taking into
account that $\|U^n\|=\|S_{nt}-A_0\|\konv 0$, as it was shown in the proof of Theorem
\ref{t2}.  
\qed

\bigskip
\textbf{Remark 2}
By means of the functions $m_x(t), M_x(t)$, which have been introduced
during the proof of Theorem \ref{t1} an estimate of the value
$\|S_t-A_0\|$ might be obtained.  

In order to show this we remember some constants (s. properties a), g)): 
$\gamma$ - the constant of nonflatness of the cone $K$ and $C_u$ a
constant which satisfies $\|x\|\le C_u\|x\|_u$ for each $x\in X$ (s. inequality (\ref{f05}).   
Finally denote by $Q_+$ the closure of the set $S_\tau(B_+)$, where  
$B_+$ is the intersection of the unit ball $B(0;1)$ with the cone $K$. 
Notice that for each $y\in K$ and $f\in\mathcal F$ one has 
\begin{equation}\label{f24}
m_y(t)\le f(S_t(y))\le M_y(t) \quad \mbox{and} \quad 
m_y(t)\le f(A_0(y))\le M_y(t).
\end{equation}
Then any vector $x\in B(0;1)$ can be represented as 
$x=\gamma\,(x'-x'')$, where $x',x''\in B_+$. 
Therefore $t>\tau$ implies 
\begin{eqnarray*}
\|S_t-A_0\| & \leq & \sup\limits_{x\in B(0;1)}\|(S_{t-\tau}-A_0)(S_\tau(x))\| \\
             & \leq & 2\gamma\,\sup\limits_{x\in B_+}\|(S_{t-\tau}-A_0)(S_\tau)(x)\|   \\
             &   =  & 2\gamma\,\sup\limits_{y\in Q_+}\|S_{t-\tau}(y)-A_0(y)\|  \\
             & \leq & 2\gamma\,C_u\sup\limits_{y\in Q_+}\sup\limits_{f\in\mathcal
F}|f(S_{t-\tau}(y))-f(A_0(y))|.
\end{eqnarray*}
In view of the inequality (\ref{f24}) we get for each $f\in\mathcal F$
\[
|f(S_{t-\tau}(y))-f(A_0(y))|\le M_y(t-\tau)-m_y(t-\tau).
\]
Together with the previous inequality this yields the required estimate
\[
\|S_t-A_0\|\le
2\gamma\,C_u\sup_{y\in Q_+}(M_y(t-\tau)-m_y(t-\tau)).
\]
One easy proves that the functions $y\mapsto M_y(t), \ y\mapsto
m_y(t)$ are continuous, and so Dini's theorem implies that the difference  
 $M_y(t-\tau)-m_y(t-\tau)$, which for $t\to+\infty$ monotonically converges to 
 zero,  uniformly decreases to $0$ on the compact set $Q_+$. In this way we
 get another proof of statement (v) of Theorem \ref{t2}.
 
\bigskip

It is easy to see that the statements of our Theorems \ref{t1} and \ref{t2}
remain to be valid also for a "discrete" semigroup of operators, i.e. for
the sequence of iterates $\{A^n\}_{n\in \mathbb{N}}$ of some positive operator $A$.  
The proofs of the Theorems \ref{t1} and \ref{t2} might be adapted to that
case, even with some obvious simplifications. Therefore we restrict
ourselves with only the formulations. 

\begin{theo}\label{t3}
Let $(X,K,\|\cdot\|)$ be an ordered Banach space and $A\in L(X)$ a 
positive operator which satisfies the following conditions   
\begin{itemize}
\item[1)] for each vector $x\in K,\,x\neq 0$ there exists a
natural number $n_x$ such that $A^{n_x}(x)\in \Int(K)$;
\item[2)] for each vector $x\in K$ its trajectory
$\{A^n(x)\}_{n\in \mathbb N}$ is relatively compact.
\end{itemize}
Then the sequence $\{A^n\}_{n\ge 0}$ pointwise converges to some 
operator $A_0$.

If this operator is not zero, then there exist a vector $u\in \Int(K)$ 
and a functional $f_0\in K'$ such that
\begin{itemize}
\item[(i)] $A(u)=u$, $A^*(f_0)=f_0, \ f_0(u)=1$ and moreover $f_0(x)>0$ if
$x\in K,\,x\neq 0$;
\item[(ii)] $A_0=f_0\otimes u$;
\item[(iii)] for each $f\in X'$ one has
$\left(A^*\right)^n(f)\konv_{n\to+\infty} A_0^*(f)$ for each
$f\in X'$ with respect to the weak*-topology $\sigma(X',X)$;
\item[(iv)] $\lambda=1$ is a simple eigenvalue of the operators
$A$ and $A^*$.
\end{itemize}
\end{theo}

\begin{theo}\label{t4}
Let $(X,K,\|\cdot\|)$ an ordered normed space and $A\in L(X)$ a 
positive operator which satisfies the following conditions   
\begin{itemize}
\item[1)] for each vector $x\in K,\,x\neq 0$ there exists a
natural number $n_x$ such that $A^{n_x}(x)\in \Int(K)$;
\item[2)] some power of  $A$ is a compact operator;
\item[3)] $\sup\limits_{n\in \mathbb N}\|A^n\|<+\infty$.
\end{itemize}
Then the statements  
of the Theorem \ref{t3} remain true. Moreover, 
\begin{itemize}
\item[(v)] $\|A^n-A_0\|\konv_{n\to\infty}0$, i.e. the operators $A^n$ 
(and, of course, the adjoint operators $(A^n)^*$) converge to   
$A_0$ (to $A_0^*$) with respect to the norm.
\end{itemize}
\end{theo}

We remark another assertion, which turns out to be a special case of
Theorem \ref{t4}. 
\begin{theo}\label{t5}
Let $(X,K,\|\cdot\|)$ an ordered normed space and $A\in L(X)$ a 
positive operator which satisfies the following conditions   
\begin{itemize}
\item[1)] some power of the operator $A$ is strongly positive, i.e. for some 
$p\in\mathbb{N}$ one has  
$A^p(K\setminus\{0\})\subset \Int(K)$; 
\item[2)] some power of  $A$ is a compact operator;
\item[3)] $\sup\limits_{n\in \mathbb N}\|A^n\|<+\infty$.
\end{itemize}
Then all statements of the Theorem \ref{t4} remain true.
\end{theo}

This theorem indeed is a special case of Theorem \ref{t4} because its
condition 1) is stronger than condition 1) of
Theorem \ref{t4} and the other assumptions are identic.

\bigskip
At the end we shall shortly deal with two examples. Let $X$ be either the vector space $C(Q)$ of all
real continuous functions defined on the compact topological space $Q$ with the cone 
$K$ of all nonnegative functions or the vector
space $\textbf{c}$ of all converging real sequences   
with the cone\footnote{\,In both cases the cone $K$ satisfies the 
condition ${\rm int}(K)\neq \varnothing$.} $K$ consisting of all nonnegative sequences. 
The symbol ${\bf 1}$ denotes correspondingly the function identically equal to $1$ on $Q$ or the 
sequence whose components are all $1$. 

An operator $A\in L(X)$ is called a {\it  Markov operator}, if it is positive and 
$A({\bf 1})={\bf 1}$.  
We indicate two examples of a compact Markov operator in the spaces 
$C([0,1])$ and $\textbf{c}$, respectively, which satisfies the condition 1) of Theorem \ref{t4} but does
not satisfy the condition 1) of Theorem \ref{t5}.

\bigskip
\textbf{Example 1} Let be $X=C([0,1])$.  
For some arbitrary fixed number 
$\theta \in(0,1)$ denote the function $\varphi(s)=(\theta+\sqrt s)^{-1}$
and put
\[
\big( A(x) \big)(s)=\varphi(s)\left(\int_0^\theta x(t)\,dt+\int_0^{\sqrt
s}x(t)\,dt\right) \qquad (x\in C([0,1])\,).
\]

Obviously, $A$ is compact operator in $C([0,1])$ and has the properties 
\[ 
A\ge 0, \quad A({\bf 1})={\bf 1}, \quad \|A\|=1 .
\]
Therefore the operator $A$ satisfies the identic conditions 2)
and 3) of the Theorems \ref{t4} and \ref{t5}.

We show that $A$ satisfies the condition 1) of Theorem \ref{t4}, but not 
the condition 1) of Theorem \ref{t5}. 
For $x\in K, \ x\not\equiv 0$ put  
\[
p=\sup\{\tau\in [0,1] \colon x(t)=0 \ \mbox{for} \ t\in [0,\tau] \}
\]
and define $p=0$ if $x(0)\neq 0$.

It is easy  to see that  $p<1$ and that $p<\theta$ implies $\big(A(x)\big)(s)>0$ on 
$[0,1]$ because of 
 $\int_0^\theta x(t)\,dt =\int_p^\theta x(t)\,dt>0$. 

If further for some $m\in \mathbb N$ the number $p$ satisfies the inequality  
\[ 
\theta^{\frac{1}{2^{m}}}\le p<\theta^{\frac{1}{2^{m-1}}}, 
\] 
then an induction argument shows
that 
\[
\big(A^{m-1}(x)\big)(0)=0, \ \mbox{ but } \ \big(A^m(x)\big)(s)>0 \ \mbox{for} \
s\in [0,1],
\]
i.e. $A^m(x)\in \Int(K)$ ($A^0=I$). Consequently, for each  $x\in K$  there exists its individual power  
$m_x$ with $A^{m_x}(x)\in \Int(K)$, however a common power, simultaneously for all $x\in K$, does not
exist. 

\medskip
\textbf{Example 2}
Let an operator $A$ be defined by a stochastic matrix as follows 
\[
A=\left(\begin{array}{llllllll}
\frac{1}{2} & \frac{1}{2}  &        0     &  0 & \ldots  & 0  & 0 &\ldots \\[2mm]
\frac{1}{2} & \frac{1}{2^2}&\frac{1}{2^2} &  0 & \ldots  & 0  & 0 & \ldots  \\[2mm]
 \vdots     &   \vdots   &  \vdots       &    &         &    &   & \\ 
 \frac{1}{2}& \frac{1}{2^2}&\frac{1}{2^3}&  \frac{1}{2^4} & \ldots  &  \frac{1}{2^n} &   
 \frac{1}{2^n} & \ldots \\
   \vdots   &    \vdots     &  \vdots     &           &   & &  &  
\end{array}\right).
\]

$A$ defines a positive operator from  $\textbf{c}$ to $\textbf{c}$, and 
obviously is a Markov operator. For each $n\in \mathbb{N}$ the matrix $A^n$ is also
stochastic. Moreover, the first $n+k$ entries in the $k$-th row of $A^n$ are positive 
($k,n=1,2,\ldots$).  
Denote by $e_n$ the sequence  $(0,\ldots,0,1,0,\ldots)$ with $1$ at the 
$n$-th position.  
Then the first coordinate of the vector $A^n(e_{n+2})$ 
(still) turns out to be zero,   
all coordinates of $A^{n+1}(e_{n+2})$ are positive and from the second one on, they are
all equal, i.e. all coordinates of $A^{n+1}(e_{n+2})$ are separated from $0$. 
For the vector $A^{n+2}(e_{n+2})$ even all coordinates are equal and positive. 
So, we have $A^n(e_{n+2})\notin \Int(K)$ but 
$A^{n+1}(e_{n+2}), A^{n+2}(e_{n+2})\in \Int(K)$. 
On the other hand is it clear that no iterate of $A$ can
satisfy the condition $A^n(K\setminus \{0\})\subset \Int (K)$. 
 
The compactness of $A$ follows from the possibility to approximate the 
operator $A- g\otimes \bf 1$ (and so the operator $A$) by finite-rank operators, 
where $g$ is the functional generated by the sequence 
$\{\frac{1}{2^n}\}_{n\in \mathbb{N}}$. 
 
The limit distribution of the operator $A$, i.e. a vector $f_0$ 
such that $A^*(f_0)=f_0$, can be calculated as the sequence with the members 
$c_n-c_{n+1}$, where 
$c_n=2^{-\frac{(n-1)n}{2}}$
for $n=1,2,\ldots$.

\bigskip
\textbf{Remark 3} We point out one particular situation, namely the situation of
Markov operators, in which Theorem \ref{t5} allows us to obtain 
as a special case some result which has been proved in 
\cite{Fel2} (chapt.VIII, \S 7, Th.1). 

Theorem \ref{t5} holds for any Markov operator $A$ in the
space $C(Q)$, if some power of $A$ is a compact operator (in particular, 
if $A$ is weakly compact) and if $A$ satisfies the condition of regularity: 
for some $m$ the inequality  $A^m(x)>0$ holds everywhere
on $Q$ for each nonnegative function $x\in C(Q)$ that is not identical zero
(with other words,  $A^m$ is strongly positive).

In this case, the functional $f_0$ whose existence and
properties are ensured by the statements of Theorem \ref{t5},
is called {\it "stationary distribution"} (s.\cite{Fel1}) or 
{\it "limit distribution"} (s.\cite{KemSn}, chapt.IV). Notice
that under the made assumptions 
for any probability measure $\mu\in C^*(Q)$
one has not only the convergence $(A^*)^n(\mu)\konv_{n\to\infty} f_0$ in 
 variation but, according to the statement ($v$) 
of Theorem \ref{t4}, even
$\sup_{\mu}\|(A^*)^n(\mu)-f_0\|\konv_{n\to\infty} 0$, i.e. the sequence
$\{(A^*)^n\}_{n\in \mathbb{N}}$ converges to the operator $A_0^*$ with respect 
to the norm. 

\bigskip
At the end we notice another remark to Theorem \ref{t4} concerning the random walk on a compact space
(compare with \cite{Fel2}, chapt.VIII, \S7 Th.1 and \cite{Bor}, \S\S6,7). We keep the terminology, which is familiar
in the theory of stationary Markov chains and also the notations there (s. \cite{Fel1}, \cite{Fel2}). 

\medskip
\textbf{Remark 4}
Let be $Q$ a metrizable compact space, $\{P_s\}_{s\in Q}\subset C^*(Q)$ a stochastic kernel 
which satisfies the conditions 
\begin{itemize}
\item[1)] for each open set $G\subset Q$ the function $s \mapsto P_s(G)$ is continuous on $Q$; 
\item[2)] for each nonempty open set $G\subset Q$ and each point $s\in Q$ there exists a number $n=n(G,s)$ 
such that $P^{(n)}_s(G)>0$, where $P^{(n)}_s$ denotes the $n$-th iterate of the kernel $P_s$. 
\end{itemize}
Then there exists a probability measure $P_0\in C^*(Q)$ such that  
\begin{itemize}
\item[(i)] $P_0(G)>0$ for any nonempty open subset $G$; 
\item[(ii)] $\sup_{s\in Q}\|P^{(n)}_s - P_0\|\konv_{n\to\infty} 0$.
\end{itemize}

For a short proof we consider the Markov operator $A$ corresponding to the kernel $\{P_s\}_{s\in Q}$, where
$(Ax)(s)=\int_Q x(t)\,{\rm d}P_s(t)$.   
Condition 1) implies that the operator $A$ is weakly compact (s. \cite{Edw}, Th.9.4.10) and consequently
the operator $A^2$ is compact (s. \cite{Edw}, sect.9.4.5). From condition 2) one gets that the operator
$A$ meets the condition 1) of Theorem \ref{t4}. Therefore, the operator $A$ satisfies all conditions of
Theorem \ref{t4}. It remains to notice that $P_0$ is that measure which is generated by the functional $f_0$.
The statement (ii) holds because of 
\[  \sup_{s\in Q} \|P^{(n)}_S-P_0\|\leq  \|(A^*)^n-A^*_0\|\konv_{n\to\infty} 0, \]  
where $A_0$ is the limit operator from Theorem \ref{t4}. 

\vspace{1cm}


\bibliography{litalt}

\begin{thebibliography}{10}

\bibitem{Bor}
{Borovkov, A.A.}
\newblock {\em Ergodicity and Stability of Random Processes}.
\newblock (Russian). Editorial URSS, Moskwa, 1999.

\bibitem{Bour1}
{Bourbaki, N.}
\newblock {\em Espaces Vectoriels Topologiques}.
\newblock Hermann, Paris, 1955.

\bibitem{CleHei}
{Clement,P., Heijmans,H.J.A.M., Angenent,S., van Duijn,C.J., de Pagter, B.}
\newblock {\em One-Parameter Semigroups}.
\newblock North Holland. CWI Monograph 5, Amsterdam, 1987.

\bibitem{Edw}
{Edwards, R.E.}
\newblock {\em Functional Analysis: Theory and Applications}.
\newblock Holt, Rinehart and Winston, New York, 1965.

\bibitem{Fel1}
{Feller, W.}
\newblock {\em An Introduction to Probability Theory and its Applications,
  Vol.I}.
\newblock Wiley, New York, 1968.

\bibitem{Fel2}
{Feller, W.}
\newblock {\em An Introduction to Probability Theory and its Applications,
  Vol.II}.
\newblock Wiley, New York, 1971.

\bibitem{KA3}
{Kantorovich, W., Akilov, G.P.}
\newblock {\em Functional Analysis}.
\newblock Pergamon Press, Oxford, 1982.

\bibitem{KemSn}
{Kemeny, J.G., Snell, J.L.}
\newblock {\em Finite Markov Chains}.
\newblock D. van Nostrand Comp.Inc., Princeton, 1960.

\bibitem{Kra}
{Krasnosel'skij, M.A.}
\newblock {\em Positive Solutions of Operator Equations,}.
\newblock P. Noordhoff, Groningen, The Netherlands, 1964.

\bibitem{KLS}
{Krasnosel'skij, M.A., Lifshits, J.A., Sobolev, A.V.}
\newblock {\em Positive Linear Systems}.
\newblock Heldermann Verlag, Berlin, 1989.

\bibitem{nev}
{Neveu, J.}
\newblock {\em The Calculus of Probability}.
\newblock Holden - Day, Inc., San Francisco, London, Amsterdam, 1965.

\end{thebibliography}
\bibliographystyle{plain}


\vspace{2cm}

\footnotesize
{\author{Boris M.~Makarow, \\ 
 Sankt Petersburg State University, 
Faculty of Mathematics and Mechanics, 98 904 Sankt Petersburg, Russia. }\\
{\it e-mail:  bm1092@gmail.com}  \\

\author{Martin R.~Weber, \\
Technische Universit\"at Dresden,
Fakult\"at f\"ur Mathematik, Institut f\"ur Analysis, 01062 Dresden, Germany.}\\
{\it e-mail: martin.weber@tu-dresden.de}
}

\end{document}